\documentclass[12pt]{article}
\usepackage{mathrsfs}
\usepackage{amssymb}
\usepackage{amsfonts}
\usepackage{latexsym}
\usepackage{amsthm}

\usepackage{mathtext}
\usepackage[T1]{fontenc}
\usepackage[ansinew]{inputenc}

\newcommand{\er}[1]{{\rm(\ref{#1})}}
\def\lb{\label}

\theoremstyle{plain}
\newtheorem{theorem}{\bf Theorem}[section]
\newtheorem{lemma}[theorem]{\bf Lemma}

\theoremstyle{remark}

\begin{document}
\def\a{\alpha}  \def\cA{{\cal A}}     \def\bA{{\bf A}}  \def\mA{{\mathscr A}}
\def\b{\beta}   \def\cB{{\cal B}}     \def\bB{{\bf B}}  \def\mB{{\mathscr B}}
\def\g{\gamma}  \def\cC{{\cal C}}     \def\bC{{\bf C}}  \def\mC{{\mathscr C}}
\def\G{\Gamma}  \def\cD{{\cal D}}     \def\bD{{\bf D}}  \def\mD{{\mathscr D}}
\def\d{\delta}  \def\cE{{\cal E}}     \def\bE{{\bf E}}  \def\mE{{\mathscr E}}
\def\D{\Delta}  \def\cF{{\cal F}}     \def\bF{{\bf F}}  \def\mF{{\mathscr F}}
\def\c{\chi}    \def\cG{{\cal G}}     \def\bG{{\bf G}}  \def\mG{{\mathscr G}}
\def\z{\zeta}   \def\cH{{\cal H}}     \def\bH{{\bf H}}  \def\mH{{\mathscr H}}
\def\e{\eta}    \def\cI{{\cal I}}     \def\bI{{\bf I}}  \def\mI{{\mathscr I}}
\def\p{\psi}    \def\cJ{{\cal J}}     \def\bJ{{\bf J}}  \def\mJ{{\mathscr J}}
\def\vT{\Theta} \def\cK{{\cal K}}     \def\bK{{\bf K}}  \def\mK{{\mathscr K}}
\def\k{\kappa}  \def\cL{{\cal L}}     \def\bL{{\bf L}}  \def\mL{{\mathscr L}}
\def\l{\lambda} \def\cM{{\cal M}}     \def\bM{{\bf M}}  \def\mM{{\mathscr M}}
\def\L{\Lambda} \def\cN{{\cal N}}     \def\bN{{\bf N}}  \def\mN{{\mathscr N}}
\def\m{\mu}     \def\cO{{\cal O}}     \def\bO{{\bf O}}  \def\mO{{\mathscr O}}
\def\n{\nu}     \def\cP{{\cal P}}     \def\bP{{\bf P}}  \def\mP{{\mathscr P}}
\def\r{\rho}    \def\cQ{{\cal Q}}     \def\bQ{{\bf Q}}  \def\mQ{{\mathscr Q}}
\def\s{\sigma}  \def\cR{{\cal R}}     \def\bR{{\bf R}}  \def\mR{{\mathscr R}}
\def\S{\Sigma}  \def\cS{{\cal S}}     \def\bS{{\bf S}}  \def\mS{{\mathscr S}}
\def\t{\tau}    \def\cT{{\cal T}}     \def\bT{{\bf T}}  \def\mT{{\mathscr T}}
\def\f{\phi}    \def\cU{{\cal U}}     \def\bU{{\bf U}}  \def\mU{{\mathscr U}}
\def\F{\Phi}    \def\cV{{\cal V}}     \def\bV{{\bf V}}  \def\mV{{\mathscr V}}
\def\P{\Psi}    \def\cW{{\cal W}}     \def\bW{{\bf W}}  \def\mW{{\mathscr W}}
\def\o{\omega}  \def\cX{{\cal X}}     \def\bX{{\bf X}}  \def\mX{{\mathscr X}}
\def\x{\xi}     \def\cY{{\cal Y}}     \def\bY{{\bf Y}}  \def\mY{{\mathscr Y}}
\def\X{\Xi}     \def\cZ{{\cal Z}}     \def\bZ{{\bf Z}}  \def\mZ{{\mathscr Z}}
\def\O{\Omega}
\def\ve{\varepsilon}
\def\vt{\vartheta}
\def\vp{\varphi}
\def\vk{\varkappa}
\def\ti{\tilde}
\def\h{\hat}

\def\Z{{\Bbb Z}}
\def\R{{\Bbb R}}
\def\C{{\Bbb C}}
\def\T{{\Bbb T}}
\def\N{{\Bbb N}}

\def\Z{{\Bbb Z}}
\def\R{{\Bbb R}}
\def\C{{\Bbb C}}
\def\T{{\Bbb T}}
\def\N{{\Bbb N}}
\def\dD{{\Bbb D}}

\def\qqq{\qquad}
\def\qq{\quad}
\def\ma{\left(\begin{array}{cc}}    \def\am{\end{array}\right)}
\def\iint{\int\!\!\!\int}
\def\lt{\biggl}                     \def\rt{\biggr}
\let\ge\geqslant                   \let\le\leqslant
\def\[{\begin{equation}}            \def\]{\end{equation}}
\def\wt{\widetilde}                 \def\pa{\partial}
\def\sm{\setminus}                  \def\es{\emptyset}
\def\no{\noindent}                  \def\ol{\overline}
\def\iy{\infty}                     \def\ev{\equiv}
\def\/{\over}
\def\we{\wedge}
\def\ts{\times}
\def\os{\oplus}
\def\su{\subset}
\def\h{\hat}
\def\wh{\widehat}
\def\Ra{\Rightarrow}
\def\ra{\rightarrow}
\def\la{\leftarrow}
\def\da{\downarrow}
\def\ua{\uparrow}
\def\lra{\leftrightarrow}
\def\Lra{\Leftrightarrow}
\def\Re{\mathop{\rm Re}\nolimits}
\def\Im{\mathop{\rm Im}\nolimits}
\def\supp{\mathop{\rm supp}\nolimits}
\def\sign{\mathop{\rm sign}\nolimits}
\def\Ran{\mathop{\rm Ran}\nolimits}
\def\Ker{\mathop{\rm Ker}\nolimits}
\def\Tr{\mathop{\rm Tr}\nolimits}
\def\const{\mathop{\rm const}\nolimits}
\def\Wr{\mathop{\rm Wr}\nolimits}

\def\th{\theta}
\def\dlint{\displaystyle\int\limits}
\def\iintt{\mathop{\int\!\!\int\!\!\dots\!\!\int}\limits}
\def\intt{\mathop{\int\int}\limits}
\def\lim{\mathop{\rm lim}\limits}
\def\mult{\!\cdot\!}
\def\BBox{\hspace{1mm}\vrule height6pt width5.5pt depth0pt \hspace{6pt}}
\def\1{1\!\!1}
\newcommand{\bwt}[1]{{\mathop{#1}\limits^{{}_{\,\bf{\sim}}}}\vphantom{#1}}
\newcommand{\bhat}[1]{{\mathop{#1}\limits^{{}_{\,\bf{\wedge}}}}\vphantom{#1}}
\newcommand{\bcheck}[1]{{\mathop{#1}\limits^{{}_{\,\bf{\vee}}}}\vphantom{#1}}
\def\nh{\bhat}
\def\nc{\bcheck}
\newcommand{\oo}[1]{{\mathop{#1}\limits^{\,\circ}}\vphantom{#1}}
\newcommand{\po}[1]{{\mathop{#1}\limits^{\phantom{\circ}}}\vphantom{#1}}
\def\ctg{\mathop{\rm ctg}\nolimits}
\def\notto{\to\!\!\!\!\!\!\!/\,\,\,}

\def\pgbrk{\pagebreak}

%%%%%%%%%%%%%%%%%%%%%%%%%%5

\def\Twelve{
\font\Tenmsa=msam10 scaled 1200 \font\Sevenmsa=msam7 scaled 1200
\font\Fivemsa=msam5 scaled 1200
%\newfam\msafam
\textfont\msafam=\Tenmsa \scriptfont\msafam=\Sevenmsa
\scriptscriptfont\msafam=\Fivemsa

\font\Tenmsb=msbm10 scaled 1200 \font\Sevenmsb=msbm7 scaled 1200
\font\Fivemsb=msbm5 scaled 1200
%\newfam\msafam
\textfont\msbfam=\Tenmsb \scriptfont\msbfam=\Sevenmsb
\scriptscriptfont\msbfam=\Fivemsb

\font\Teneufm=eufm10 scaled 1200 \font\Seveneufm=eufm7 scaled 1200
\font\Fiveeufm=eufm5 scaled 1200
%\newfam\eufmfam
\textfont\eufmfam=\Teneufm \scriptfont\eufmfam=\Seveneufm
\scriptscriptfont\eufmfam=\Fiveeufm}

\def\Ten{
\textfont\msafam=\tenmsa \scriptfont\msafam=\sevenmsa
\scriptscriptfont\msafam=\fivemsa

\textfont\msbfam=\tenmsb \scriptfont\msbfam=\sevenmsb
\scriptscriptfont\msbfam=\fivemsb

\textfont\eufmfam=\teneufm \scriptfont\eufmfam=\seveneufm
\scriptscriptfont\eufmfam=\fiveeufm}

\title{The parametrization of the Marchenko-Ostrovsky mapping in terms of
the Dirichlet eigenvalues}

\author{Maria Evgenievna Korotyaeva\begin{footnote}
{Mathematisches Institut, Humboldt-Universit\"{a}t zu Berlin,
e-mail: korotiaeva@yahoo.de}\end{footnote}}

\maketitle

\begin{abstract}
\no We consider the inverse spectral problem for periodic Jacobi
matrices in terms of the vertical slits on the quasi-momentum domain
plus the Dirichlet eigenvalues, i.e., the Marchenko-Ostrovsky
mapping. Moreover, we show that the gradients of the Dirichlet
eigenvalues and of the so-called norming constants are linear
independent.
\end{abstract}

\vskip 0.25cm
\section {Introduction}
\setcounter{equation}{0}

We consider the self-adjoint N-periodic Jacobi operator $\cJ$ on a
Hilbert space $\ell^2=\ell^2(\Z)$ given by $(\cJ y)_{n\in\Z}=
a_{n}y_{n+1}+b_{n}y_{n}+a_{n-1}y_{n-1}$ for
$y=(y_n)_{n\in\Z}\in\ell^2$ and for the $N$-periodic sequences
$a_{n}=e^{x_n}>0,\ x_n,b_{n}\in\R.$ Furthermore, we assume
\[\label{1.4}
q=(x,b)\equiv(x_n,b_n)_{n\in\Z_N}\in\mH^2,\ \mH\equiv\{b\in\R^N:
\sum_{n=1}^Nb_n=0\}.
\]
For simplicity, we use the notation $\Z_N=\{1,2,...,N\},\ N\in\N,$
throughout this paper. To begin, we recall some well known facts
(see, e.g., \cite{vM}). Let $\vp=\vp(\l,q)=(\vp_n(\l,q))_{n\in\Z}$
and $\vt=\vt(\l,q)=(\vt_n(\l,q))_{n\in\Z}$ denote two fundamental
solutions of the equation
\[\label{1.5}
a_{n-1}y_{n-1}+b_ny_n+a_ny_{n+1}=\l y_{n},\ (\l,n)\in\C\times\Z,
\]
with the initial conditions $\vp_0\equiv\vt_1\equiv 0,\
\vp_1\equiv\vt_0\equiv 1.$ The Lyapunov function
$\D(\l,q)={1\/2}(\vp_{N+1}(\l,q)+\vt_N(\l,q))$ is the discriminant
of the equation \er{1.5} and characterizes the spectrum
$\s(q)=\{\l\in\R: |\D(\l,q)|\leq 1\}$ of $\cJ.$ The spectrum of
$\cJ$ is absolutely continuous and consists of $N$ bands
$\s_n=[\l_{n-1}^{+}, \l_{n}^{-}],\ n\in\Z_N,$ separated by the gaps
$\g_n=(\l_n^{-},\l_n^{+}),$ where $\l_n^{\pm}=\l_n^{\pm}(q)$ are the
roots of $\D^2(\l,q)=1$ and satisfy
$\l_N^{+}\equiv\l_0^{+}<\l_1^{-}\leq\l_1^{+}<...<\l_{N-1}^{-}\leq\l_{N-1}^{+}<\l_N^{-}.$
That is, if a gap $\g_n$ degenerates, then the corresponding
segments $\s_n,\ \s_{n+1}$ merge. Moreover, there is exactly one
point $\l_n=\l_n(q)\in[\l^{-}_n,\l^{+}_n]$ for each $n\in\Z_{N-1}$
such that
\[\label{1.9}
\D^{'}(\l_n,q)=0,\ \D^{''}(\l_n,q)\neq 0,\
(-1)^{N-n}\D(\l_n,q)\geq1.
\]
Here and below, $(')=\partial/\partial\l$.

Traditionally since \cite{vM}, the inverse spectral problem for the
periodic Jacobi operator has been solved using the Neumann
eigenvalues $\m_n=\m_n(q),\ n\in\Z_{N-1},$ given by the zeroes of
$\vt_{N+1}(\l,q)=0.$ The main goal of this paper is to solve the
inverse spectral problem alternatively using the Dirichlet
eigenvalues $\n_n=\n_n(q),\ n\in\Z_{N-1},$ given by the zeroes of
$\vp_{N}(\l,q)=0.$ That is, we define the auxiliary spectrum by the
Dirichlet spectrum instead of the Neumann spectrum. Note that
$\m_n,\n_n\in [\l^{-}_n,\l^{+}_n], n\in\Z_{N-1}.$

To outline the plan of this note, we recall that the inverse
spectral problem consists of the following four parts, namely,

\no i) \textit{The uniqueness.} Prove that the spectral data
uniquely determines the potential.

\no ii) \textit{The reconstruction.} Give an algorithm to recover
the potential from the spectral data.

\no iii) \textit{The characterization.} Give the conditions for some
data to be the spectral data of some potential.

\no iv) \textit{The stability.} Give the two-sided a priori
estimates of the potential in terms of the spectral data.

\no We construct a Marchenko-Ostrovsky mapping
$\p:\mH^2\ra\R^{2N-2}$ for the periodic Jacobi operator in terms of
the Dirichlet eigenvalues $\n_n$ by $\p=(\p_n)_{n\in\Z_{N-1}},\
\p_n=(\p_{1,n},\p_{2,n})\in\R^2,$ where
\[\lb{9.11}
\p_{1,n}=\log((-1)^{N-n}\vp_{N+1}(\n_n)),\
\p_{2,n}=(|\p_n|^2-\p_{1,n}^2)^{1\/2}\sign(\l_n-\n_n),
\]
\[\lb{9.12}
\cosh|\p_n|=(-1)^{N-n}\D(\l_n).
\]
Here, $\p_{1,n}$ is the so-called norming constant. It is easy to
verify $|\p_n|^2-|\p_{1,n}|^2\geq 0$ since \er{9.11}, \er{9.12} and
the Wronskian identity $\vp_{N+1}\vt_N-\vp_N\vt_{N+1}=1$ together
imply
\[\lb{9.13}
(-1)^{N-n}\D(\n_n,q)=\cosh\p_{1,n}(q).
\]
Note incidentally that the mapping $\p$ is an analogue of the
Marchenko-Ostrovski mapping \cite{MO} for the continuous case and
has similar properties (see \cite{MO}, \cite{Ko}).

Firstly, we will prove the characterization and the uniqueness
showing that the mapping $\{$potential$\}$ $\ra$ $\{$spectral
data$\}$ is a homeomorphism.

\begin{theorem}
\label{T.97} The mapping $\p:\mH^2\ra\R^{2N-2}$ is a real analytic
isomorphism between the Hilbert spaces $\mH^2$ and $\R^{2N-2}$.
\end{theorem}

\no {\bf Remark.} We recall some necessary definitions. Let $\mH,
\mH_0$ be Hilbert spaces. The derivative of a map $f:\mH\to \mH_0$
at a point $q\in \mH$ is a bounded linear map from $\mH$ into
$\mH_0$, which we denote by $d_qf$. A map $f:\mH\to\mH_0$ is a real
analytic isomorphism between $\mH$ and $\mH_0$ if $f$ is bijective
and both $f$ and $f^{-1}$ are real analytic maps of the space.

Secondly, we will obtain the reconstruction and the stability. We
use for it the geometric interpretation of the Marchenko-Ostrovsky
mapping, which is similar to the continues case mentioned in
\cite{MO} and \cite{Ko}. For this purpose, we introduce the
conformal mapping (the quasi-momentum) $\kappa:\L\to K$ by
\[\lb{ca}
\cos\kappa(\l)=(-1)^N\D(\l,q),\ \ \ \l\in\L,\ \ \ \ \ {\rm and}\ \ \
\ \ \ \kappa(it)\to\pm i\iy\ \ \ {\rm as}\ \ t\to\pm\iy.
\]
Here and below, $\L=\C\sm\cup_1^{N-1}\g_n$ is the domain,
$K=\{\kappa:0\le\Re\kappa\le N\pi\}\sm\cup_1^{N-1}\kappa(\g_n)$ is
called the quasi-momentum domain and $\G_n=(\pi n+i|\p_n|,\pi
n-i|\p_n|)$ is an excised vertical slit.

\begin{theorem}
\label{T.2} i) For each $\p\in\R^{2N-2},$ there is a unique point
$q\in\mH^2$ and a unique conformal mapping $\kappa:\L\to K$ such
that the following identities hold true
\[\lb{ci}
\kappa(\l_n(q)\pm i0)=\pi n\pm i|\p_n(q)|,\ \ \ \kappa(\n_n(q)\pm
i0)=\pi n\pm i\p_{1,n}(q),\ \ \ n\in\Z_{N-1}.
\]
ii) For $(\p_{1,n})_{n\in\Z_{N-1}},$ there is a standard algorithm
to recover $a,b.$

\no iii) The following two-sided estimates hold true
\[\lb{ciii}
{1\/4}e^{2\max|\p_n|}<{1\/4}(\l_N^{-}-\l_0^{+})^2<b^2+2a^2<N(\l_N^{-}-\l_0^{+})^2
<16Ne^{2\max |\p_n|}.
\]
\end{theorem}

A crucial argument in the proof of Theorem \ref{T.97} is the fact
that the gradients of the Dirichlet eigenvalues and of the norming
constants are linear independent. More precisely, we define the
symplectic form $\wedge$ by
\[\label{9.4}
f\wedge g=\sum_{n=1}^N(f_{1,n}g_{2,n}-f_{2,n}g_{1,n})
-(f_{1,n-1}g_{2,n}-f_{2,n}g_{1,n-1})
\]
for $f=(f_{1,n},f_{2,n})_{n\in\Z_N},\
g=(g_{1,n},g_{2,n})_{n\in\Z_N}\in\C^{2N}$ with $f_{1,0}\equiv
f_{1,N}$ and $g_{1,0}\equiv g_{1,N}.$ Note that below $\d_{n,m}$
stands for the Kronecker symbol for all $n,m\in\Z$. Then we show

\begin{theorem}
\label{T.95} For all $n,m\in\Z_{N-1},$ it holds true
\[\lb{95.1}
d_q\n_n\wedge d_q\n_m=0,
\]
\[\lb{95.2}
d_q\p_{1,n}\wedge d_q\p_{1,m}=0,
\]
\[\lb{95.3}
d_q\p_{1,n}\wedge d_q\n_m=2\d_{n,m}.
\]
In particular, $d_q\n_n,\ d_q\p_{1,n},\ n\in\Z_{N-1},$ is a basis of
$\mH^2.$
\end{theorem}

\no P\"{o}schel and Trubowitz \cite{PT} proved an analogue of
Theorem \ref{T.95} for the Sturm-Liouville problem on the interval
$[0,1].$ We use their arguments in our proof. Note that van Moerbeke
\cite{vM} proved (using another approach) that the gradients of the
Neumann eigenvalues $\m_n$ and of the norming constants are linear
independent. Namely, van Moerbeke used the Jacobi matrices with
removed rows and columns. Remark that our proof also can be applied
for the case of $\m_n.$

There are different approaches to the inverse spectral problem for
the periodic Jacobi operator. The investigation on this topic
started in 1976 by van Moerbeke \cite{vM} and by Date and Tanaka
\cite{DT}. Both works obtained the reconstruction, but not the
characterization: van Moerbeke did it using the Stieltjes inverse
spectral method from \cite{Ah} or \cite{GK}, and Date and Tanaka did
it applying the suffix shifting by a constant. The nonlinear Toda
lattice turned out an important application of these methods (see
\cite{To}). The first for us known work on the characterization is
the paper \cite{Pe} by Perkolab, where some analogue of Theorem
\ref{T.2} is showed using \cite{MO}. Further, Korotyaev and Kutsenko
\cite{KoKu} showed Theorems \ref{T.97} and \ref{T.2} in terms of the
Neumann eigenvalues applying approach \cite{KaKo}. That is,
\cite{KoKu} extended the result of Marchenko and Ostrovski about the
height-slit mapping for the Hill operator (see \cite{MO},
\cite{Ko1}) to the case of the periodic Jacobi matrix using the
Neumann eigenvalues. Lastly, the inverse problem in terms of the gap
lengths was solved in \cite{BGGK} based on the approach from
\cite{GT} and in \cite{Ko1} based on the approach from \cite{KaKo}.

Our note is organized as follows: Section 2 displays some
preliminary statements in terms of $\n_n.$ In Section 3 we prove
Theorem \ref{T.95}; this is technically the most difficult part of
this note. In Section 4, we show Theorem \ref{T.97} using the
argument from \cite{Ko} and \cite{KoKu}, where this theorem in terms
of $\m_n$ is proved. Then Theorem \ref{T.97} together with
\cite{KoKr} and \cite{KoKu} directly implies Theorem \ref{T.2}.

\section{Preliminaries}
\setcounter{equation}{0}

In this section, we will determine the gradients (with respect to
$q$) of the Dirichlet eigenvalues and the norming constants in terms
of the fundamental solutions. For this purpose, we define the
Wronskian by
\[\label{9.5}
\{f,g\}_n=a_n(f_ng_{n+1}-f_{n+1}g_n),\ n\in\Z,
\]
for the sequences $f=(f_n)_{n\in\Z},\ g=(g_n)_{n\in\Z}$ with $f_n,
g_n\in\C.$ Below, we use the notation
$\partial=\partial_q=(\partial_{x_k},\partial_{b_k})_{k\in\Z_N}.$

\begin{lemma}
\label{L.90} Each from the functions $\n_n,\ \p_{1,n},\
n\in\Z_{N-1},$ is real analytic on $\mH^2$ and satisfies
\[\lb{9.14}
d_q\n_n=-{\partial\vp_{N}(\n_n(q),q)\/\vp'_{N}(\n_n(q),q)},
\]
\[\lb{9.17}
d_q\p_{1,n}={\vp'_{N+1}(\n_n(q),q)d_q\n_n+\partial\vp_{N+1}(\n_n(q),q)\/
\vp_{N+1}(\n_n(q),q)}.
\]
\end{lemma}
\no {\bf Proof.} This proof is similar to the continues case
\cite{Ko} (see also \cite{KoKu}). \BBox

\begin{lemma}
\label{L.91} Let $\h{\vp}=\vp(\n_n(q),q), \h{\vt}=\vt(\n_n(q),q)$
for all $(n,q)\in\Z_{N-1}\times\mH^2.$ Then for all $k\in\Z_N$, the
following identities hold
\[\lb{9.1}
d_{q_k}\n_n=-{(2a_k\h{\vp}_k\h{\vp}_{k+1},\h{\vp}_k^2)\/
a_N\h{\vp}_{N+1}\h{\vp}^{'}_{N}},
\]
\[\lb{9.2}
d_{q_k}\p_{1,n}=B_{n,k}+(\h{\vp}^{'}_{N+1}\h{\vt}_{N}
-\h{\vp}^{'}_{N}\h{\vt}_{N+1})d_q\n_n,
\]
\[\lb{9.3}
B_{n,k}={1\/a_N}(a_k(\h{\vp}_{k+1}\h{\vt_k}
+\h{\vp}_k\h{\vt}_{k+1}),\h{\vp}_k\h{\vt}_k).
\]
\end{lemma}

\no {\bf Proof.} We assume $k,j\in\Z_N$.

i) We want to show \er{9.1} applying \er{9.14}. That is, we have to
determine the derivation $\partial_{q_k}\h{\vp}_{N}.$ Firstly, we
calculate $\partial_{x_k}\h{\vp}_{N}$ using the equation \er{1.5}
for $\vp_j$
$$
a_{j-1}\vp_{j-1}+(b_j-\l)\vp_j+a_j\vp_{j+1}=0,
$$
and its derivation with respect to $x_k$
$$
a_{j-1}\partial_{x_k}\vp_{j-1}
+(b_j-\l)\partial_{x_k}\vp_j+a_j\partial_{x_k}\vp_{j+1}
=-a_k(\delta_{j,k}\vp_{k+1}+\delta_{j,k+1}\vp_k+\delta_{k,N}\delta_{j,1}\vp_0).
$$
Below, $\chi_{k<N}$ stands for the characteristic function, i.e.
$\chi_{k<N}=1$ for $k<N$ and $\chi_{k<N}=0$ for $k\geq N$ (recalling
that $N$ is fixed). Multiplying the first equation by
$\partial_{x_k}\vp_j$ and the second one by $\vp_j$ and taking the
difference, we sum the result over all $j\in\Z_N,$ that is,
$$
2a_k\vp_k\vp_{k+1}=a_k(\chi_{k<N}2\vp_k\vp_{k+1}+\delta_{k,N}(\vp_0\vp_1+\vp_N\vp_{N+1}))
$$$$
=\sum_{j=1}^Na_{j-1}(\vp_{j-1}\partial_{x_k}\vp_j-\vp_j\partial_{x_k}\vp_{j-1})
+a_j(\vp_{j+1}\partial_{x_k}\vp_j-\vp_{j}\partial_{x_k}\vp_{j+1})
$$$$
=\sum_{j=1}^N\{\partial_{x_k}\vp,\vp\}_j-\{\partial_{x_k}\vp,\vp\}_{j-1}
=\{\partial_{x_k}\vp,\vp\}_N-\{\partial_{x_k}\vp,\vp\}_{0}
=\{\partial_{x_k}\vp,\vp\}_N
$$
since $\vp_0=0.$ Next, setting $\l=\n_n$ and recalling
$\h{\vp}_{N}=0,$ it gives
\[\lb{91.1}
\partial_{x_k}\h{\vp}_{N}={2a_k\h{\vp}_k\h{\vp}_{k+1}\/a_N\h{\vp}_{N+1}}.
\]
Secondly, we calculate $\partial_{b_k}\h{\vp}_{N}$ using the
equation \er{1.5} for $\vp_j$
$$
a_{j-1}\vp_{j-1}+(b_j-\l)\vp_j+a_j\vp_{j+1}=0,
$$
and its derivation with respect to $b_k$
$$
a_{j-1}\partial_{b_k}\vp_{j-1}
+(b_j-\l)\partial_{b_k}\vp_j+a_j\partial_{b_k}\vp_{j+1}
=-\delta_{j,k}\vp_{k}.
$$
Multiplying the first equation by $\partial_{b_k}\vp_j$ and the
second one by $\vp_j$ and taking the difference, we sum the result
over all $j\in\Z_N,$ that is,
$$
\vp_k^2=\sum_{j=1}^Na_{j-1}(\vp_{j-1}\partial_{b_k}\vp_j-\vp_j\partial_{b_k}\vp_{j-1})
+a_j(\vp_{j+1}\partial_{b_k}\vp_j-\vp_j\partial_{b_k}\vp_{j+1})
$$$$
=\sum_{j=1}^N\{\partial_{b_k}\vp,\vp\}_j-\{\partial_{b_k}\vp,\vp\}_{j-1}
=\{\partial_{b_k}\vp,\vp\}_N-\{\partial_{b_k}\vp,\vp\}_{0}
=\{\partial_{b_k}\vp,\vp\}_N
$$
since $\vp_0=0.$ Setting
$\l=\n_n,$ we see that $\h{\vp}_{N}=0$ implies
$\partial_{b_k}\h{\vp}_{N}={\h{\vp}_k^2\/a_N\h{\vp}_{N+1}}$ and by
\er{91.1}, we get \er{9.1}.

ii) In order to prove \er{9.2} and \er{9.3}, we determine
$\partial_{q_k}\h{\vp}_{N+1}$ and then substitute it into the
identity \er{9.17}. Firstly, we calculate
$\partial_{x_k}\h{\vp}_{N+1}$ using the equation \er{1.5} for
$\vp_j$
$$
a_{j-1}\vp_{j-1}+(b_j-\l)\vp_j+a_j\vp_{j+1}=0,
$$
and the derivation of the equation \er{1.5} for $\vt_j$ with respect
to $x_k$
$$
a_{j-1}\partial_{x_k}\vt_{j-1}+(b_j-\l)\partial_{x_k}\vt_j+a_j\partial_{x_k}\vt_{j+1}
=-a_k(\delta_{j,k}\vt_{k+1}+\delta_{j,k+1}\vt_k+\delta_{k,N}\delta_{j,1}\vt_0).
$$
Multiplying the first equation by $\partial_{x_k}\vt_j$ and the
second one by $\vp_j$ and taking the difference, we sum the result
over all $j\in\Z_N,$ that is,
$$
a_k(\chi_{k<N}(\vp_k\vt_{k+1}+\vp_{k+1}\vt_k)+\delta_{k,N}(\vp_1\vt_{0}+\vp_{N+1}\vt_N))
$$$$
=\sum_{j=1}^Na_{j-1}(\vp_{j-1}\partial_{x_k}\vt_j-\vp_j\partial_{x_k}\vt_{j-1})
+a_j(\vp_{j+1}\partial_{x_k}\vt_j-\vp_j\partial_{x_k}\vt_{j+1})
$$$$
=\sum_{j=1}^N\{\partial_{x_k}\vt,\vp\}_j-\{\partial_{x_k}\vt,\vp\}_{j-1}
=\{\partial_{x_k}\vt,\vp\}_N-\{\partial_{x_k}\vt,\vp\}_{0}
=\{\partial_{x_k}\vt,\vp\}_N
$$
since $\partial_{x_k}\vt_0=\partial_{x_k}\vt_1=0.$ Setting
$\l=\n_n,$ we get
$$
a_k(\h\vp_k\h\vt_{k+1}+\h\vp_{k+1}\h\vt_k)=\{\partial_{x_k}\h\vt,\h\vp\}_N.
$$
We observe that $\h{\vp}_{N}=0$ implies
$$
a_k(\h{\vp}_k\h{\vt}_{k+1}+\h{\vp}_{k+1}\h{\vt}_k)
=a_N\h{\vp}_{N+1}\partial_{x_k}\h{\vt}_{N}
=-a_N\h{\vt}_{N}\partial_{x_k}\h{\vp}_{N+1}
+a_N\h{\vt}_{N+1}\partial_{x_k}\h{\vp}_{N}
$$
since $\partial_{x_k}\{\vt,\vp\}_N=0$ and hence
\[\lb{91.2}
\partial_{x_k}\h{\vp}_{N+1}
=-{a_k\/a_N}(\h{\vt}_N(\h{\vt}_k\h{\vp}_{k+1}+\h{\vp}_k\h{\vt}_{k+1})
+2\h{\vt}_{N+1}\h{\vp}_k\h{\vp}_{k+1}).
\]
Secondly, we calculate $\partial_{b_k}\h{\vp}_{N+1}$ using the
equation \er{1.5} for $\vp_j$
$$
a_{j-1}\vp_{j-1}+(b_j-\l)\vp_j+a_j\vp_{j+1}=0,
$$
and the derivation of the equation \er{1.5} for $\vt_j$ with respect
to $b_k$
$$
a_{j-1}\partial_{b_k}\vt_{j-1}
+(b_j-\l)\partial_{b_k}\vt_j+a_j\partial_{b_k}\vt_{j+1}
=-\delta_{j,k}\vt_{k}.
$$
Multiplying the first equation by $\partial_{b_k}\vt_j$ and the
second one by $\vp_j$ and taking the difference, we sum the result
over all $j\in\Z_N,$ that is,
$$\vp_{k}\vt_k=
\sum_{j=1}^Na_{j-1}(\vp_{j-1}\partial_{b_k}\vt_j-\vp_j\partial_{b_k}\vt_{j-1})
+a_j(\vp_{j+1}\partial_{b_k}\vt_j-\vp_j\partial_{b_k}\vt_{j+1})
$$$$
=\sum_{j=1}^N\{\partial_{b_k}\vt,\vp\}_j-\{\partial_{b_k}\vt,\vp\}_{j-1}
=\{\partial_{b_k}\vt,\vp\}_N-\{\partial_{b_k}\vt,\vp\}_{0}
=\{\partial_{b_k}\vt,\vp\}_N
$$ since
$\partial_{b_k}\vt_0=\partial_{b_k}\vt_1=0.$ We set $\l=\n_n,$ then
$\h{\vp}_{N}=0$ implies
$$
\h{\vp}_{k}\h{\vt}_k=a_N\h{\vp}_{N+1}\partial_{b_k}\h{\vt}_{N}
=-a_N\h{\vt}_{N}\partial_{b_k}\h{\vp}_{N+1}
+a_N\h{\vt}_{N+1}\partial_{b_k}\h{\vp}_{N},
$$
by $\partial_{b_k}\{\vt,\vp\}_N=0.$ Therefore,
$\partial_{b_k}\h{\vp}_{N+1}
=a_N^{-1}(\h{\vp}_{N+1}\h{\vp}_k\h{\vt}_k
-\h{\vt}_{N+1}\h{\vp}_k^2)$ and \er{91.2} yield
$$
\partial_{q_k}\h{\vp}_{N+1}=
{1\/a_N}\lt(\h{\vp}_{N+1}(a_k(\h{\vp}_{k+1}\h{\vt}_k
+\h{\vp}_k\h{\vt}_{k+1}),\h{\vp}_k\h{\vt}_k)
-\h{\vt}_{N+1}(2a_k\h{\vp}_k\h{\vp}_{k+1},\h{\vp}_k^2)\rt).
$$
Substituting this expression into the identity \er{9.17}, one
obtains \er{9.2} and \er{9.3}. \BBox

\section{The proof of Theorem \ref{T.95}.}
\setcounter{equation}{0}

In this Section, we will show that $d_q\n_n,\ d_q\p_{1,n},\
n\in\Z_{N-1},$ is a basis of $\mH^2$ using the symplectic form
defined in \er{9.4}. First of all, we display some identities for
the Wronskian defined in \er{9.5}.

\begin{lemma}
\label{L.92} Let $\h{\vp}=\vp(\n_n(q),q), \h{\vt}=\vt(\n_n(q),q)$
and $\tilde{\vp}=\vp(\n_m(q),q), \tilde{\vt}=\vt(\n_m(q),q),\
n,m\in\Z_{N-1}, q\in\mH^2.$ Then for all $k\in\Z_N,$ the following
identities hold true
\[\lb{92.1}
\sum_{j=1}^N\h{\vp}_j\tilde{\vp}_j=0,\ \
\{\tilde{\vp},\h{\vp}\}_k=(\n_n-\n_m)\sum_{i=1}^k\h{\vp}_i\tilde{\vp}_i,
\]
\[\lb{92.2}
\sum_{j=1}^N\h{\vt}_j\tilde{\vt}_j={\{\tilde{\vt},\h{\vt}\}_N\/(\n_n-\n_m)},\
\
\{\tilde{\vt},\h{\vt}\}_k=(\n_n-\n_m)\sum_{i=1}^k\h{\vt}_i\tilde{\vt}_i,
\]
\[\lb{92.3}
\sum_{j=1}^N\h{\vt}_j\tilde{\vp}_j={a_0(1-\tilde{\vp}_{N+1}\h{\vt}_N)\/(\n_n-\n_m)},\
\
\{\tilde{\vp},\h{\vt}\}_k=-a_0+(\n_n-\n_m)\sum_{i=1}^k\h{\vt}_i\tilde{\vp}_i,
\]
\[\lb{92.4}
\sum_{j=1}^N\h{\vp}_j\tilde{\vt}_j={a_0(\tilde{\vt}_N\h{\vp}_{N+1}-1)\/(\n_n-\n_m)},\
\
\{\tilde{\vt},\h{\vp}\}_k=a_0+(\n_n-\n_m)\sum_{i=1}^k\h{\vp}_i\tilde{\vt}_i,
\]
\[\lb{92.5}
\sum_{j=1}^N\h{\vp}_j^2=a_0\h{\vp}_{N+1}\h{\vp}_{N}^{'}.
\]
\end{lemma}

\no {\bf Proof.} We assume $n,m\in\Z_{N-1}$ and $j,k\in\Z_N$.

i) To show \er{92.1}, we consider the equation \er{1.5} for
$\h{\vp}_j$ and $\tilde{\vp}_j$ respectively
$$
a_{j-1}\h{\vp}_{j-1}+b_j\h{\vp}_j+a_j\h{\vp}_{j+1}=\n_n\h{\vp}_{j},
$$$$
a_{j-1}\tilde{\vp}_{j-1}+b_j\tilde{\vp}_j+a_j\tilde{\vp}_{j+1}=\n_m\tilde{\vp}_j.
$$
Multiplying the first equation by $\tilde{\vp}_j$ and the second one
by $\h{\vp}_j$ and taking the difference, we get
$$
a_{j-1}(\tilde{\vp}_j\h{\vp}_{j-1}-\tilde{\vp}_{j-1}\h{\vp}_j)
+a_j(\tilde{\vp}_j\h{\vp}_{j+1}-\tilde{\vp}_{j+1}\h{\vp}_j)=(\n_n-\n_m)\h{\vp}_j\tilde{\vp}_j
$$
or equivalently, by the definition \er{9.5} of the Wronskian,
\[\lb{9.6}
\{\tilde{\vp},\h{\vp}\}_j-\{\tilde{\vp},\h{\vp}\}_{j-1}
=(\n_n-\n_m)\h{\vp}_j\tilde{\vp}_j.
\]
Summing \er{9.6} over all $j\in\Z_N$ and using
$\tilde{\vp}_0=\h{\vp}_0=\tilde{\vp}_{N}=\h{\vp}_{N}=0,$ we get the
first identity in \er{92.1}
$$
(\n_n-\n_m)\sum_{j=1}^N\h{\vp}_j\tilde{\vp}_j=
\{\tilde{\vp},\h{\vp}\}_N-\{\tilde{\vp},\h{\vp}\}_0=0.
$$
Next, we show the second identity in \er{92.1} by induction on $k.$
For $k=1,$ the identity holds according to \er{9.6} since
$\{\tilde{\vp},\h{\vp}\}_0=0.$ Assuming that the identity holds for
$k-1,$ we will verify it for $k:$ By \er{9.6} and the induction
hypothesis, it follows
$$
{\{\tilde{\vp},\h{\vp}\}_k\/(\n_n-\n_m)}={\{\tilde{\vp},\h{\vp}\}_{k-1}\/(\n_n-\n_m)}
+\h{\vp}_k\tilde{\vp}_k=\sum_{i=1}^{k-1}\h{\vp}_i\tilde{\vp}_i
+\h{\vp}_k\tilde{\vp}_k=\sum_{i=1}^k\h{\vp}_i\tilde{\vp}_i.
$$
This establishes the second identity in \er{92.1}.

ii) The proof of \er{92.2} is similar to that of \er{92.1} since
$\{\tilde{\vt},\h{\vt}\}_0=0.$

iii) To show \er{92.3}, we consider again the equation \er{1.5} for
$\h{\vt}_j$ and $\tilde{\vp}_j$ respectively
$$
a_{j-1}\h{\vt}_{j-1}+b_j\h{\vt}_j+a_j\h{\vt}_{j+1}=\n_n\h{\vt}_{j},
$$$$
a_{j-1}\tilde{\vp}_{j-1}+b_j\tilde{\vp}_j+a_j\tilde{\vp}_{j+1}=\n_m\tilde{\vp}_{j}.
$$
Multiplying the first equation by $\tilde{\vp}_j$ and the second one
by $\h{\vt}_j$ and taking the difference, we get
$$
a_{j-1}(\tilde{\vp}_j\h{\vt}_{j-1}-\tilde{\vp}_{j-1}\h{\vt}_j)
+a_j(\tilde{\vp}_j\h{\vt}_{j+1}-\tilde{\vp}_{j+1}\h{\vt}_j)=(\n_n-\n_m)\h{\vt}_j\tilde{\vp}_j
$$
or equivalently
\[\lb{9.7}
\{\tilde{\vp},\h{\vt}\}_j-\{\tilde{\vp},\h{\vt}\}_{j-1}
=(\n_n-\n_m)\h{\vt}_j\tilde{\vp}_j.
\]
Next, summing \er{9.7} over all $j\in\Z_N$ and using
$\{\tilde{\vp},\h{\vt}\}_0=-a_0,\ a_0=a_N,\ \h{\vp}_{N}=0,$ the
first statement in \er{92.2} follows
$$
(\n_n-\n_m)\sum_{j=1}^N\h{\vt}_j\tilde{\vp}_j=
\{\tilde{\vp},\h{\vt}\}_N-\{\tilde{\vp},\h{\vt}\}_0=a_0(1-\tilde{\vp}_{N+1}\h{\vt}_N).
$$
We prove the second identity \er{92.2} by induction on $k.$ For
$k=1,$ it holds according to \er{9.7}, since
$\{\tilde{\vp},\h{\vt}\}_0=-a_0$ and $\tilde{\vp}_1\h{\vt}_1=0.$
Supposing the identity is truth for $k-1,$ prove it for $k.$
\er{9.7} and the induction hypothesis together imply
$$
\{\tilde{\vp},\h{\vt}\}_k=
\{\tilde{\vp},\h{\vt}\}_{k-1}+(\n_n-\n_m)\h{\vt}_k\tilde{\vp}_k=
(\n_n-\n_m)\sum_{i=1}^{k-1}\h{\vt}_i\tilde{\vp}_i+(\n_n-\n_m)\h{\vt}_k\tilde{\vp}_k-a_0=
$$$$
=(\n_n-\n_m)\sum_{i=1}^k\h{\vt}_i\tilde{\vp}_i-a_0.
$$
So the second identity in \er{92.2} follows.

iv) The proof of \er{92.4} is similar to that of \er{92.3}, since we
have $\{\tilde{\vt},\h{\vp}\}_0=a_0$ and
$\{\tilde{\vt},\h{\vp}\}_N=a_N\tilde{\vt}_N\h{\vp}_{N+1}=a_0\tilde{\vt}_N\h{\vp}_{N+1}.$

v) To verify \er{92.5}, we use the equation \er{1.5} for $\vp_j(\l)$
and its gradient with respect to $\l$
$$
a_{j-1}\vp_{j-1}(\l)+(b_j-\l)\vp_j(\l)+a_j\vp_{j+1}(\l)=0,
$$$$
a_{j-1}\vp^{'}_{j-1}(\l)+(b_j-\l)\vp^{'}_j(\l)+a_j\vp^{'}_{j+1}(\l)=\vp_{j}(\l).
$$
Multiplying the first equation by $\vp^{'}_j(\l)$ and the second one
by $\vp_j(\l)$ and taking the difference, we get
$$
a_{j-1}(\vp^{'}_j(\l)\vp_{j-1}(\l)-\vp^{'}_{j-1}(\l)\vp_j(\l))
+a_j(\vp^{'}_j(\l)\vp_{j+1}(\l)-\vp^{'}_{j+1}(\l)\vp_j(\l))=\vp^2_j(\l)
$$
or equivalently
\[\lb{9.8}
\{\vp^{'}(\l),\vp(\l)\}_j-\{\vp^{'}(\l),\vp(\l)\}_{j-1}
=\vp^2_j(\l).
\]
Next, summing \er{9.8} over all $j\in\Z_N$ and using
$\vp_1\equiv\vp^{'}_1\equiv 0,$ we obtain
$$
\sum_{j=1}^N\vp^2_j(\l)=
\{\vp^{'}(\l),\vp(\l)\}_N-\{\vp^{'}(\l),\vp(\l)\}_0
=\{\vp^{'}(\l),\vp(\l)\}_N.
$$
Setting $\l=\n_n$ and recalling $\h{\vp}_{N}=0,$ it becomes
$\sum_{j=1}^N\h{\vp}^2_j=\{\h{\vp}^{'},\h{\vp}\}_N
=a_N\h{\vp}_N\h{\vp}^{'}_{N}.$ \BBox

\ \

Now we will apply the definition \er{9.4} of the symplectic form:
For $n\in\Z_{N-1}$, we consider $d_q\n_n=(d_{q_k}\n_n)_{k\in\Z_N}$
given in \er{9.1} and $B_n=(B_{n,k})_{k\in\Z_N}$ defined in
\er{9.3}. We observe that $d_{q_0}\n_n\equiv d_{q_N}\n_n$ and
$B_{n,0}\equiv B_{n,N}$ for all $n\in\Z_{N-1}$ since $q_0=q_N.$ So
we can define a symplectic form for $d_q\n_n$ and $B_n$.

\begin{theorem}
\label{T.93} For all $n,m\in\Z_{N-1},$ the following identities hold
\[\lb{93.1}
d_q\n_n\wedge d_q\n_m=0,
\]
\[\lb{93.2}
B_n\wedge B_m=0,
\]
\[\lb{93.3}
B_n\wedge d_q\n_m=-2\d_{n,m}.
\]
\end{theorem}

\no {\bf Proof.} We assume $n,m\in\Z_{N-1}$. Furthermore, we use the
following abbreviations $\h\vp=\vp(\n_n(q),q), \h\vt=\vt(\n_n(q),q)$
and $\ti{\vp}=\vp(\n_m(q),q), \ti{\vt}=\vt(\n_m(q),q)$ for all
$q\in\mH^2.$

i) To verify the identity \er{93.1}, let $n\neq m$ since the case
$n=m$ is obvious. Applying the definition \er{9.4} of the symplectic
form and introducing
$C=2(a_N^2\h\vp_{N+1}\h\vp_{N}^{'}\tilde{\vp}_{N+1}\tilde{\vp}_{N}^{'})^{-1},$
we have
$$
d_q\n_n\wedge d_q\n_m
=C\sum_{k=1}^N\lt(a_{k}(\h\vp_{k}\h\vp_{k+1}\tilde{\vp}_{k}^2
-\h\vp_{k}^2\tilde{\vp}_{k}\tilde{\vp}_{k+1})
-a_{k-1}(\h\vp_{k-1}\h\vp_{k}\tilde{\vp}_{k}^2
-\h\vp^2_{k}\tilde{\vp}_{k-1}\tilde{\vp}_{k})\rt).
$$
Using the definition \er{9.5} of the Wronskian, this equals
$$
d_q\n_n\wedge d_q\n_m =C\sum_{k=1}^N\h\vp_{k}\tilde{\vp}_{k}
(\{\tilde{\vp},\h\vp\}_{k}+\{\tilde{\vp},\h\vp\}_{k-1}).
$$
By \er{92.1}, it follows
$$
{d_q\n_n\wedge d_q\n_m\/C(\n_n-\n_m)}
=\sum_{k=1}^N\h\vp_{k}\tilde{\vp}_{k}
(\sum_{i=1}^{k-1}\h\vp_i\tilde{\vp}_i+\sum_{i=1}^k\h\vp_i\tilde{\vp}_i)
=\sum_{k=1}^N\h\vp_{k}\tilde{\vp}_{k}\sum_{i=1}^{k-1}\h\vp_i\tilde{\vp}_i
+\sum_{k=1}^N\h\vp_{k}\tilde{\vp}_{k}\sum_{i=1}^k\h\vp_i\tilde{\vp}_i.
$$
Below, we need the following simple identities
\[\lb{R1}
\sum_{k=1}^Nz_k\sum_{i=1}^kw_i=\sum_{k=1}^Nw_k\sum_{i=k}^Nz_i,\ \ \
\ \ \
\sum_{k=2}^Nz_k\sum_{i=1}^{k-1}w_i=\sum_{k=1}^{N-1}w_k\sum_{i=k+1}^Nz_i
\]
for all $z=(z_k)_{k\in\Z_N}, w=(w_k)_{k\in\Z_N}\in\C^N.$ Then we can
use \er{R1} and \er{92.5} to obtain
$${d_q\n_n\wedge d_q\n_m\/C(\n_n-\n_m)}
=\sum_{k=1}^N\h\vp_k\tilde{\vp}_{k}\sum_{i=1}^{k-1}\h\vp_i\tilde{\vp}_i
+\sum_{k=1}^N\h\vp_{k}\tilde{\vp}_{k}\sum_{i=k}^N\h\vp_i\tilde{\vp}_i
=\sum_{k=1}^N\h\vp_{k}\tilde{\vp}_{k}\sum_{i=1}^{N}\h\vp_i\tilde{\vp}_i=0,
$$
which vanishes by the sum identity in \er{92.1}.

ii) We verify the identity \er{93.2} in the same manner like
\er{93.1}. Again, we consider $n\neq m,$ since the case $n=m$ is
obvious, and get
$$
B_n\wedge B_m={1\/a_N^2}\sum_{k=1}^N
a_k\lt((\h\vp_{k+1}\h\vt_k+\h\vp_k\h\vt_{k+1})\tilde{\vp}_k\tilde{\vt}_k
-\h\vp_k\h\vt_k(\tilde{\vp}_{k+1}\tilde{\vt}_k+\tilde{\vp}_k\tilde{\vt}_{k+1})\rt)
$$$$
-a_{k-1}\lt((\h\vp_k\h\vt_{k-1}+\h\vp_{k-1}\h\vt_k)\tilde{\vp}_k\tilde{\vt}_k
-\h\vp_k\h\vt_k(\tilde{\vp}_k\tilde{\vt}_{k-1}+\tilde{\vp}_{k-1}\tilde{\vt}_k)\rt)
$$$$
={1\/a_N^2}\sum_{k=1}^N a_k\lt(
\tilde{\vp}_k\tilde{\vp}_k(\h\vt_{k+1}\tilde{\vt}_k-\h\vt_k\tilde{\vt}_{k+1})
+\h\vt_k\tilde{\vt}_k(\h\vp_{k+1}\tilde{\vp}_k-\h\vp_k\tilde{\vp}_{k+1})\rt)
$$$$
+a_{k-1}\lt(
\h\vt_k\tilde{\vt}_k(\h\vp_k\tilde{\vp}_{k-1}-\h\vp_{k-1}\tilde{\vp}_k)
+\h\vp_k\tilde{\vp}_k(\h\vt_k\tilde{\vt}_{k-1}-\h\vt_{k-1}\tilde{\vt}_k)\rt)
$$$$
={1\/a_N^2}(\sum_{k=1}^N
\h\vt_{k}\tilde{\vt}_{k}(\{\tilde{\vp},\h\vp\}_{k}+\{\tilde{\vp},\h\vp\}_{k-1})
+\sum_{n=1}^N
\h\vp_{k}\tilde{\vp}_{k}(\{\tilde{\vt},\h\vt\}_{k}+\{\tilde{\vt},\h\vt\}_{k-1})).
$$
Using firstly representations of the Wronskian in \er{92.1} and
\er{92.2} and applying \ref{R1}, we obtain
$$
{a_N^2(B_n\wedge B_m)\/(\n_n-\n_m)}
=\sum_{k=1}^N\h\vt_{k}\tilde{\vt}_{k}
(\sum_{i=1}^{k-1}\h\vp_i\tilde{\vp}_i+\sum_{i=1}^{k}\h\vp_i\tilde{\vp}_i)
+\sum_{k=1}^N\h\vp_{k}\tilde{\vp}_{k}
(\sum_{i=1}^{k-1}\h\vt_i\tilde{\vt}_i+\sum_{i=1}^{k}\h\vt_i\tilde{\vt}_i)
$$$$
=\sum_{k=1}^N\h\vt_{k}\tilde{\vt}_{k}\sum_{i=1}^{N}\h\vp_i\tilde{\vp}_i
+\sum_{k=1}^N\h\vp_k\tilde{\vp}_k\sum_{i=1}^{k-1}\h\vt_i\tilde{\vt}_i=0,
$$
which vanishes by the sum identity in \er{92.1}.

iii) The proof of the identity \er{93.3} is similar. We have
$$
B_n\wedge d_q\n_m=
(a_N^2\tilde{\vp}_{N+1}\tilde{\vp}_{N}^{'})^{-1}\sum_{k=1}^N
a_k((\h\vt_k\h\vp_{k+1}+\h\vt_{k+1}\h\vp_k)\tilde{\vp}_k^2
-2\h\vp_k\h\vt_k\tilde{\vp}_k\tilde{\vp}_{k+1})
$$$$
-a_{k-1}((\h\vt_{k-1}\h\vp_k+\h\vt_k\h\vp_{k-1})\tilde{\vp}_k^2
-2\h\vp_k\h\vt_k\tilde{\vp}_{k-1}\tilde{\vp}_k).
$$
If $n=m,$ then $\{\h\vt,\h\vp\}_k=a_N, k\in\Z,$ and \er{92.5} imply
$$
B_n\wedge d_q\n_n=(a_N^2\h\vp_{N+1}\h\vp_{N}^{'})^{-1}\sum_{k=1}^N
a_{k-1}(\h\vt_{k-1}\h\vp_k-\h\vt_k\h\vp_{k-1})\vp_k^2
+a_k(\h\vt_k\h\vp_{k+1}-\h\vt_{k+1}\h\vp_k)\h\vp_k^2
$$$$
=2(a_N^2\h\vp_{N+1}\h\vp_{N}^{'})^{-1}\sum_{k=1}^Na_N\h\vp_k^2=2.
$$
If $n\neq m,$ then for
$C=(a_N^2\tilde{\vp}_{N+1}\tilde{\vp}_{N}^{'})^{-1},$ we get
$$
B_n\wedge d_q\n_m=C\sum_{k=1}^N
a_k(\h\vp_k\tilde{\vp}_k(\tilde{\vp}_k\h\vt_{k+1}-\tilde{\vp}_{k+1}\h\vt_k)
+\h\vt_k\tilde{\vp}_k(\tilde{\vp}_k\h\vp_{k+1}-\tilde{\vp}_{k+1}\h\vp_k))
$$$$
+a_{k-1}(\h\vp_k\tilde{\vp}_k(\tilde{\vp}_{k-1}\h\vt_k-\tilde{\vp}_k\h\vt_{k-1})
+\h\vt_k\tilde{\vp}_k(\tilde{\vp}_{k-1}\h\vp_k-\tilde{\vp}_k\h\vp_{k-1}))
$$$$
=C\lt(\sum_{k=1}^N\h\vp_{k}\tilde{\vp}_{k}
(\{\tilde{\vp},\h\vt\}_{k-1}+\{\tilde{\vp},\h\vt\}_{k})
+\sum_{k=1}^N\h\vt_{k}\tilde{\vp}_{k}
(\{\tilde{\vp},\h\vp\}_{k-1}+\{\tilde{\vp},\h\vp\}_{k})\rt).
$$
Using firstly the representations of the Wronskians in \er{92.1} and
\er{92.4} and $-2a_0\sum_{k=1}^N\h\vp_k\tilde{\vp}_k=0$ due to
\er{92.1}, this gives by \er{R1}
$$
B_n\wedge d_q\n_m=(\n_n-\n_m)C\lt( \sum_{k=1}^N\h\vp_k\tilde{\vp}_k
(\sum_{i=1}^{k-1}\h\vt_i\tilde{\vp}_i+\sum_{i=1}^k\h\vt_i\tilde{\vp}_i)
+\sum_{k=1}^N\h\vt_k\tilde{\vp}_k
(\sum_{i=1}^{k-1}\h\vp_i\tilde{\vp}_i+\sum_{i=1}^k\h\vp_i\tilde{\vp}_i)
\rt)
$$$$
=(\n_n-\n_m)C\lt(
\sum_{k=1}^N\h\vp_k\tilde{\vp}_k\sum_{i=1}^{N}\h\vt_i\tilde{\vp}_i
+\sum_{k=1}^N\h\vp_k\tilde{\vp}_k\sum_{i=1}^{N}\h\vt_i\tilde{\vp}_i\rt)=0,
$$
which vanishes by the sum identity in \er{92.1}. \BBox

The last Theorem allows us to prove the main result of this Section:

\no {\bf Proof of Theorem \ref{T.95}.} Firstly, Theorem \ref{T.93}
implies the required equations: \er{95.1} is the statement
\er{93.1}. Next, \er{95.2} follows by
$$
d_q\p_{1,n}\wedge d_q\p_{1,m}=
(B_n+(\h\vp^{'}_{N+1}\h\vt_{N}-\h\vp^{'}_{N}\h\vt_{N+1})d_q\n_n)\wedge
(B_m+(\tilde{\vp}^{'}_{N+1}\tilde{\vt}_{N}-\tilde{\vp}^{'}_{N}\tilde{\vt}_{N+1})d_q\n_m)
$$$$
=(\h\vp^{'}_{N+1}\h\vt_{N}-\h\vp^{'}_{N}\h\vt_{N+1})(2\d_{n,m}-2\d_{n,m})=0.
$$
Applying \er{93.1} and \er{93.3}, we obtain
$$
d_q\p_{1,n}\wedge d_q\n_m=
(B_n+(\h\vp^{'}_{N+1}\h\vt_{N}-\h\vp^{'}_{N}\h\vt_{N+1})d_q\n_n)\wedge
d_q\n_m=B_n\wedge d_q\n_m=2\d_{n,m}.
$$
Secondly, the identities \er{95.1}--\er{95.3} yield that $d_q\n_n,
d_q\p_{1,n},\ n\in\Z_{N-1},$ is a basis of $\mH^2.$ \BBox

\section{The Proof of Theorems \ref{T.97} and \ref{T.2}}
\setcounter{equation}{0}

In order to prove Theorem \ref{T.97}, we need the following lemma:

\begin{lemma}
\label{L.9} Each from the functions $\l_n,\ \xi_n\equiv|\p_n|^2,\
\p_n,\ n\in\Z_{N-1},$ is real analytic on $\mH^2$ and satisfies
\[\lb{9.15}
d_q\l_n=-{\partial\D^{'}(\l_n(q),q)\/\D^{''}(\l_n(q),q)},
\]
\[\lb{9.16}
d_q\xi_n=(-1)^{N-n}{\partial\D(\l_n(q),q)\/2(d\cosh\sqrt{\xi_n}/d\xi_n)},
\]
\[\lb{9.18}
(-1)^{N-n}(\sinh
\p_{1,n})d_q\p_{1,n}=\partial\D(\n_n(q),q)+\D'(\n_n(q),q)d_q\n_n.
\]
Moreover, there exists a real analytic positive function $\b_n$ on
$\mH^2$ such that for all $q\in\mH^2$,
\[\lb{9.19}
\p_{2,n}(q)=\b_n(q)(\l_n(q)-\n_n(q)).
\]
\end{lemma}
\no {\bf Proof.} This proof is similar to the continues case
\cite{Ko} (see also \cite{KoKu}). \BBox

\no {\bf Proof of Theorem \ref{T.97}.} Although we use the approach
from the papers \cite{KaKo} and \cite{Ko}, we give the accurate
proof for the sake of the reader. We need the following theorem from
the nonlinear functional analysis; it is a modification \cite{KoKu}
of a "basic theorem" of the direct method in \cite{KaKo}.

\begin{theorem}
\label{T.96} Let $\mH,$ and $\mH_0$ be Hilbert spaces equipped with
norms $\|\cdot\|$ and $\|\cdot\|_0$ respectively. Let
$f_0:\mH\ra\mH_0$ be a real analytic isomorphism between $\mH$ and
$\mH_0.$ If a map $f:\mH\ra\mH_0$ satisfies following conditions:

i) $f$ is a local real analytic isomorphism,

ii) $f-f_0$ is a compact map, i.e., it maps a weakly convergent
sequence in $\mH$ into a convergent sequence in $\mH_0,$

iii) We have $\|f(x)\|_0\ra\iy$ as $\|x\|\ra\iy$ and it holds
$f^{-1}(0)=0.$

\no Then $f$ is a real analytic isomorphism between $\mH$ and
$\mH_0.$
\end{theorem}

\no We will check the conditions of Theorem \ref{T.96} for the
Marchenko-Ostrovsky mapping $\p$ using lemma \ref{L.9}:

i) Let us verify the first condition applying the Inverse Function
Theorem. By Lemma \ref{L.9}, $\p(\cdot)$ is real analytic on
$\mH^2$. It remains to show that $d_q\p$ is invertible. We will
prove it by contradiction. Let a vector $v\in\mH^2$ be the solution
of the equation
$$
(d_q\p)v=0\ \Lra\ \{\langle d_q\p_{1,n},v\rangle =0,\ \langle
d_q\p_{2,n},v\rangle =0,\ n\in\Z_{N-1}\}
$$
for some $q\in\mH^2$. Here, $\langle q,\widetilde{q}\rangle
=\sum_{n=1}^N(a_n\widetilde{a_n}+b_n\widetilde{b_n})$ denotes the
inner product in $\R^{2N}.$ The function
$\xi_n=\p_{1,n}^2+\p_{2,n}^2$ is analytic and the definition of $v$
implies $\langle d_q\xi_n,v\rangle=0,\ n\in\Z_{N-1}.$ Define the
polynomial $f(\l)\equiv\langle(\partial\D)(\l,q),v\rangle,\
\l\in\C,$ of degree $N-1$ with respect to $\l.$ Then \er{9.16}
implies
$$
f(\l_n)=(-1)^{N-n}{d\cosh\sqrt\xi_n\/d\xi_n}\langle
d_q\xi_n,v\rangle=0
$$
for all $n\in\Z_{N-1}.$ Therefore, $f\equiv 0.$ For fixed
$q\in\mH^2$ there are three cases:

\no 1) Let $\p_{2,n}=0.$ The differentiation of the equation
\er{9.19} implies $d_s\p_{2,n}=\b_n(s)(d_s\l_n(s)-d_s\n_n(s)),$ if
$\p_{2,n}(s)=0$ for some $s\in\mH^2.$ Using the definition of $v$
and \er{9.15} and $f=0,$ we obtain $\langle d_q\l_n,v\rangle=0$ and
hereby $\langle d_q\n_n,v\rangle=0$ for all $n\in\Z_{N-1}.$

\no 2) Let $\p_{1,n}\neq 0,\ \p_{2,n}\neq 0.$ Then \er{9.11} and
\er{9.12} yield $\l_n\neq\n_n.$ Moreover, we may apply \er{9.18} and
$f=0$ to get
$$
0=(-1)^{N-n}\sinh\p_{1,n}\langle d_q\p_{1,n},v\rangle
=\D'(\n_n,q)\langle d_q\n_n,v\rangle,
$$
that is, $\langle d_q\n_n,v\rangle=0$ since $\D'(\n_n,q)\neq 0.$

\no 3) Let $\p_{1,n}=0\neq\p_{2,n}.$ By \er{9.18}, we have
$\partial\D(\n_n,q)=-\D'(\n_n,q)d_q\n_n.$ The equation \er{9.13}
implies $\vt_{N}(\n_n,q)=\vp_{N+1}(\n_n,q)=(-1)^{N-n},$ which gives
$\l_n\neq \n_n,$ that is, $\D'(\n_n(q),q)\neq 0.$ I.e. $\langle
d_q\n_n,v\rangle=0$ since $f\equiv 0.$ The vectors
$\{d_q\p_{1,n},d_q\n_n\}_{n\in\Z_{N-1}}$ are a basis of $\mH^2$
according to Theorem \ref{T.95}, then $v=0$ holds and the operator
$d_q\p$ is invertible.

ii) The second condition, namely, the compactness follows since
$\mH$ and $\mH_0$ are finitely dimensional.

iii) The third condition follows from \cite{KoKu}, p. 6-7, where the
analog statement in terms of $\m_n$ is proved, since the norm of
$\p$ for $\m_n$ and $\n_n$ is the same. In particular, the stability
\er{ciii} is essential.

Now we see that all the conditions of Theorem \ref{T.96} are
fulfilled, then $\p$ is a real analytic isomorphism between $\mH^2$
and $\R^{2N-2}.$ \BBox

\ \

\no {\bf Proof of Theorem \ref{T.2}.} i) By Theorem \ref{T.97}, for
each $\p\in\R^{2N-2}$ there exists exactly one point $q\in\mH^2$
such that \er{9.11}-\er{9.13} hold. Next, for any point $q\in\mH^2$
there is exactly one conformal mapping $\kappa:\L\to K $ with the
properties \er{ca} (see \cite{KoKr}), which together with
\er{9.11}-\er{9.13} gives \er{ci}.

ii) The function $\kappa:\L\to K$ satisfies the equations
$\kappa(\s_n)=[\pi (n-1),\pi n], n\in\Z_N$ and $\kappa(\g_n)=\G_n,
n\in\Z_{N-1}.$ That is, if we know $(|\p_n|)_1^{N-1}$, then we get
$\L$. Moreover, \er{ci} gives all $\n_n, n\in\Z_{N-1}.$

iii) The claim \er{ciii} is proved in terms of $\m_n$ in
\cite{KoKu}, p. 2, based on \cite{KoKr}. This proof holds for the
case of $\n_n$ since $|\p_n|$ is independent of $\n_n$ or $\m_n$
according to its definition in \er{9.12}. \BBox

%\no {\bf Acknowledgements.}

\end{document}